\documentstyle[12pt]{article}
\date{}
\newtheorem{th}{Theorem}
\newtheorem{prop}{Proposition}
\newtheorem{cor}{Corollary}
\newtheorem{lem}{Lemma}

\newtheorem{dfn}{Definition}

\makeatletter

\renewcommand{\@begintheorem}[2]{\begin{trivlist}\it
 \item[\hspace{\labelsep}{\bf #1\ #2.}]}
\renewcommand{\@endtheorem}{\end{trivlist}}
\makeatother

\title{The generalized Chern character and
 \\ Lefschetz numbers in W*-modules}
\author{A.A. Pavlov
\thanks{ The work is partially supported by the
RFBR grant 99-01-01202 and INTAS grant 96-1099.}}

\begin{document}
\maketitle

\def\lra{\longrightarrow}
\def\im{\mathop{\rm Im}\nolimits}
\def\Ker{\mathop{\rm Ker}\nolimits}
\def\pcom{\;\Box}

\vspace{-1cm}
 \begin{center}
\large March 20, 2000
 \end{center}
\vspace{0.1cm}

{\small
\begin{center}
{\bf Abstract }
 \end{center}

We define $N$-theory being some analogue of $K$-theory on the category
of von Neumann algebras such that $K_0(A)\subset N_0(A)$ for
any von Neumann algebra $A$. Moreover, it turns out to be possible to
construct the extension of the Chern character to some homomorphism
from $N_0(A)$ to even Banach cyclic homology of $A$. Also, we define
generalized Lefschetz numbers for an arbitrary unitary endomorphism $U$ of
an $A$-elliptic complex. We study them in the situation when
$U$ is an element of a representation of some compact Lie group.


}

\section{Introduction}
For an arbitrary von Neumann algebra $A$ we introduce an abelian
group $N_0(A)$ in the following way. It is possible to define
some equivalence relation between normal elements of the
inductive limit $M_\infty (A)=\lim\limits_{\lra} M_n (A)$
such that for projections it coinside with the usual stable
equivalence relation.
Then the set of all equivalence classes of normal elements from
$M_\infty (A)$ is an abelian semigroup (with respect to the direct
sum operation) and $N_0(A)$ is its symmetrization.
The first part of our paper is devoted to the consideration
of some propeties of $N$-groups.
More detail on this subject can be found in~\cite{Pav2}.

Further, we introduce Banach cyclic homology of $A$
as some analogue of usual cyclic homology and
construct the generalized Chern character as a map from
$N_0(A)$ to even Banach cyclic homology.
Furthermore, this map is an extension of the classic Chern
character to the group $N_0(A)\supset K_0(A)$ in some natural sence.

In the final section we define generalized Lefschetz numbers
for an arbitrary unitary endomorphism $U$ of an $A$-elliptic
complex. Besides, in the case when $U$ is an element of
a representation of some compact Lie group we describe the connection
of the generalized Lefschetz numbers with
the W*-Lefschetz numbers of the first and of the second types
introduced in~\cite{Tro,Tro1,Fra-Tro}.

Some results of the present paper were formulated in~\cite{Pav1}.\\

The paper is organized as follows:

\begin{tabbing}
5. \= The generalized Chern character \kill
1. \> Introduction\\
2. \> Some properties of $N$-groups\\
3. \> The group $N_0(A)_{fin}$\\
4. \> Banach cyclic homology\\
5. \> The generalized Chern character\\
6. \> Generalized Lefschetz numbers\\
   \>References
 \end{tabbing}

\section{Some properties of $N$-groups}
Suppose $A$ is a von Neumann algebra,
$M_r(A)$ is the set of $r\times r$ matrices with entries in $A$,
$M_{\infty}(A)$ is the inductive limit of the sequence
$\{M_r(A)\}_{r=1}^\infty$,
and $M_{\infty}(A)_{\nu}$ is the set of normal elements for
$M_{\infty}(A)$. Denote by ${\cal B}({\bf C})$ the family
of all Borel subsets of the complex plane.
If $a\in M_{\infty}(A)_{\nu}$ and $E\in {\cal B}({\bf C})$,
then by $P_a(E)$ we denote the spectral projection of $a$
corresponding to the set $E$.
We remark that $P_a(E)\in M_{\infty}(A)_{\nu}$ since
von Neumann algebras are closed with respect to the Borel
calculus.
We denote the stable equivalence relation (see \cite{Mur})
of projections $p,q\in M_{\infty}(A)$ by $p\simeq q$.
Finally, let $sp(a)$ denote the spectrum of an element $a$.

A Borel set $E\subset{\bf C}$ is called {\em admissible\/} if
zero does not belong to the closure of $E$. Denote by
${\cal B}_*({\bf C})$ the family of all admissible Borel subsets
of the complex plane.

\begin{dfn}\label{dfn1}
 {\rm Call elements $a,b\in M_\infty(A)$}  equivalent
 {\rm (and denote by $a\simeq b$)
if and only if $P_a(E)\simeq P_b(E)$
for all $E\in {\cal B}_*({\bf C})$.}
\end{dfn}

Note that this equivalence relation coinsides with the usual
stable equivalence relation whenever $a,b$ are projections.
It is easy to see that $a\oplus 0_m\simeq a$, where
$0_m$ is the zero $m\times m$ matrix and $a\in M_{\infty}(A)_{\nu}$.
We put
$$
{\cal N}(A)=M_{\infty}(A)_{\nu} /{\simeq}.
$$
For $a\in M_{\infty}(A)_{\nu}$ let us denote
by $[a]$ the equivalence class of $a$ in ${\cal N}(A)$.
Since
$$
P_{a\oplus b}(E)=P_a(E)\oplus P_b(E)
$$
for all $a,b\in M_{\infty}(A)_{\nu}, E\in {\cal B}({\bf C})$,
this implies that the set ${\cal N}(A)$ is an abelian semigroup
with respect to the direct sum operation.

\begin{dfn}
{\rm The symmetrization of ${\cal N}(A)$ is called the}
$N$-group {\rm of $A$ and is denoted by  $N_0(A)$. }
\end{dfn}

Under the previous considerations, the following result is clear.

\begin{prop}
$K_0(A)$ is a subgroup of the group $N_0(A)$.
$\pcom$
\end{prop}

 \begin{prop}
If the group $K_0(A)$ is trivial, then the group $N_0(A)$ is trivial too.
 \end{prop}

{\bf Proof.}
Suppose $[a]-[b]\in N_0(A)$. Since $K_0(A)$ is trivial, this implies that
$P_a(E)\simeq 0\simeq P_b(E)$ for all $E\in {\cal B}_*({\bf C})$.
Whence, $a\simeq b$. Thus the group $N_0(A)$ is trivial.
$\pcom$

Note that ${\cal N}(A)$ is a {\em cancellation semigroup,\/}
i.e., the condition $[a]+[c]=[b]+[c]$ implies $[a]=[b]$ for any
$[a],[b],[c]\in {\cal N}(A)$.
In particular, the symmetrization homomorphism
$s:{\cal N}(A)\lra N_0(A),\; s([a])=[a]-[0]$
is injective.

Our next aim is to establish a functorial property for $N_0$.
We recall that an arbitrary *-homomorphism of C*-algebras
is a contraction, i.e., its norm does not
exeed $1$~\cite[Theorem 2.1.7]{Mur}. Besides,
an arbitrary surjective *-homomorphism of von Neumann algebras
is continuous with respect to the ultra-strong
topology~\cite[Theorem 2.4.23]{BR}.

Now let $A,B$ be von Neumann algebras and $\varphi: A\lra B$
an ultra-strong continuous unital *-homomorphism.
By definition, put
$\varphi (a)=(\varphi (a_{ij}))$ and
$\varphi_* ([a])=[\varphi (a)]$
for each matrix $a=(a_{ij})\in M_{\infty}(A)_{\nu}$.

\begin{th}\label{prop2}
The map $\varphi_*:N_0(A)\lra N_0(B)$ is a well defined homomorphism
of abelian groups.
\end{th}

The following lemma is the main ingredient of the proof
of Theorem~\ref{prop2}.

\begin{lem}\label{lem2}
$P_{\varphi (a)}(E)=\varphi (P_a(E))$ for each  $a\in M_{\infty}(A)_{\nu}$
and for each Borel subset $E\subset sp(a)$.
\end{lem}

{\bf Proof.}
We can assume that $a\in M_r(A)$ for some $r\ge 1$.
It is clear that $sp(\varphi (a))\subset sp(a)$.
It can be directly verified that $\varphi (R(a))=R(\varphi (a))$
for an arbitrary polynomial $R$.
Further, for any function $f$, which is continuous on $sp(a)$,
there exists a sequence of polinomials
$\{R_n\}_{n=1}^{\infty}$ such that it converges uniformly to $f$. Then
$$
 \begin{array}{ccl}
\varphi (f(a))&=&\varphi (\lim_n R_n(a))=\lim_n\varphi (R_n(a)) \\
&=&\lim_n R_n (\varphi (a))=f(\varphi (a)).
 \end{array}
$$
Now let $\chi_E$ be the characteristic function of $E$.
Then we can find a sequence $\{f_n\}_{n=1}^{\infty}$ of
continuous functions on the compact space $sp(a)$ such that
$\{f_n\}_{n=1}^{\infty}$ converges to $\chi_E$ pointwise,
i.e., with respect to the strong topology.
Moreover, we can assume that the family $\{f_n\}_{n=1}^\infty$
is norm-bounded.
In this case the sequence $\{f_n(a)\}_{n=1}^\infty$
of elements of the von Neumann algebra $M_r(A)$
converges strongly to the element $\chi_E(a)=P_a(E)$
and the family $\{f_n(a)\}_{n=1}^\infty$ is norm-bounded.
Since strong and ultra-strong topologies coinside on bounded sets,
we have the convergence with respect to the ultra-strong topology:
$P_a(E)= {\sigma}\mbox{-}\!\lim_n f_n(a)$.
Finally, using the ultra-strong continuity of the *-homomorphism
$\varphi$, we obtain:
$$
 \begin{array}{rcl}
\varphi(P_a(E))&=&\varphi({\sigma}\mbox{-}\!\lim_n f_n(a))=
{\sigma}\mbox{-}\!\lim_n\varphi (f_n (a)) \\
&=&{\sigma}\mbox{-}\!\lim_n f_n(\varphi (a))=
\chi_E (\varphi (a))=P_{\varphi (a)}(E).
\pcom
 \end{array}
$$

{\bf Proof of Theorem \ref{prop2}.}
We have to establish that $\varphi_*$ is well defined.
Let elements $a,b\in M_{\infty}(A)_{\nu}$ be equivalent. Then
$\varphi (P_a(E))\simeq\varphi (P_b(E))$
for any $E\in {\cal B}_*({\bf C})$.
Now it immediately follows from Lemma~\ref{lem2}
that the elements $\varphi (a)$ and $\varphi (b)$ are equivalent too.
$\pcom$

\section{The group $N_0(A)_{fin}$}
Let $M_{\infty}(A)_{fin}\subset M_{\infty}(A)_{\nu}$ be the subset
of all elements $a\in M_\infty (A)$ such that their spectrum is finite.
Suppose,
$$
{\cal N}(A)_{fin}=\{[a] : a\in M_{\infty}(A)_{fin}\}.
$$
We remark that ${\cal N}(A)_{fin}$ is a subsemigroup of
${\cal N}(A)$. By $N_0(A)_{fin}$ we denote the symmetrization
of the abelian monoid ${\cal N}(A)_{fin}$.
So $N_0(A)_{fin}$ is a subgroup of the group $N_0(A)$.
Any element of the group $N_0(A)_{fin}$ we can
represent in the form
$
[\bigoplus_{i=1}^n\lambda_i p_i]-[\bigoplus_{i=1}^n\lambda_i q_i],
$
where $\lambda_i\in {\bf C}$ and $p_i,q_i$ are some
(possibly zero) projections in $M_\infty (A)$.

For an arbitrary map $f:{\bf C}\setminus\{0\}\lra K_0(A)$
let us put
$\Lambda_f=\{\lambda\in {\bf C}\setminus\{0\} : f(\lambda)\ne 0\}$.
Let us denote by
$$
{\cal M}_{fin;A}={\cal M}_{fin}({\bf C}\setminus\{0\},K_0(A))
$$
the set of all maps from ${\bf C}\setminus\{0\}$ to $K_0(A)$
such that $\Lambda_f$ is finite (or empty).
The set ${\cal M}_{fin;A}$ is an abelian group
with respect to the pointwise addition of maps.
Let us consider the map
$$
\phi:N_0(A)_{fin}\lra {\cal M}_{fin;A},\quad
(\phi([a]-[b]))(\lambda)=[P_a(\{\lambda\})]-[P_b(\{\lambda\})].
$$

\begin{th}\label{th1}
The map $\phi$ is an isomorphism of groups.
\end{th}

{\bf Proof.}
It is clear that $\phi$ is a well defined homomorphism.
Let $[a]-[b]\in  N_0(A)_{fin}$, and $\phi ([a]-[b])=0$.
Therefore for each $\lambda\in {\bf C}\setminus\{0\}$
the projections $P_a(\{\lambda\})$ and $P_b(\{\lambda\})$
are stably equivalent.
Since the spectra of elements $a,b$ are finite,
we conclude that these elements are equivalent.
Hence, $\phi$ is injective.

Let us examine a map $f\in {\cal M}_{fin;A}$ such that
$\Lambda_f=\{\lambda_i\}_{i=1}^n$ and
$f(\lambda_i)=[p_i]-[q_i]$ $(1\le i\le n)$,
where $p_i,q_i$ are projections from $M_\infty (A)$.
Since $p_i\simeq 0_k\oplus p_i$, we can assume that
the projections $\{p_i\}_{i=1}^n$ (and $\{q_i\}_{i=1}^n$)
are pairwise orthogonal. We put
$a=\bigoplus_{i=1}^n \lambda_i p_i$ and
$b=\bigoplus_{i=1}^n \lambda_i q_i$.
Then the elements $a,b$ belong to $M_\infty (A)_{fin}$.
Furtermore, $\phi ([a]-[b])=f$.
Hence, $\phi$ is surjective.
$\pcom$

 \begin{cor}
The groups $N_0(M_r({\bf C}))$ and
${\cal M}_{fin}({\bf C}\setminus\{0\},{\bf Z})$ are isomorphic.
 \end{cor}

{\bf Proof.}
The spectrum of any element from $M_\infty({\bf C})$
is finite. Therefore, $N_0(M_r({\bf C}))=N_{fin}(M_r({\bf C}))$.
To complete the proof, it remains to use Theorem~\ref{th1}.
$\pcom$

Assume
$h_{\lambda ,\mu ;p}:=
[p(\lambda +\mu )]-[p\lambda\oplus p\mu ]$ and
$g_{\lambda ;p}^{(n)}:=
[\lambda p^{\oplus n}]-[n\lambda p]$,
where $n\in {\bf N}$, $\lambda ,\mu\in {\bf C}$ and
$p$ is a projection in  $M_{\infty}(A)$.
Then $h_{\lambda ,\mu ;p}$ and $g_{\lambda ;p}^{(n)}$
are elements from $N_0(A)_{fin}$.

Let us denote by $H$ the subgroup of $N_0(A)_{fin}$ with
the following system of generators
$ \{h_{\lambda ,\mu ;p} :
\lambda ,\mu\in {\bf C},\;
p \mbox{ is a projection in } M_{\infty}(A)\},
$
and by $G$ the subgroup of $N_0(A)_{fin}$ with
the following system of generators
$ \{g_{\lambda ;p}^{(n)} :
n\in{\bf N}, \lambda\in {\bf C},\;
p \mbox{ is a projection in } M_{\infty}(A)\}.
$

\begin{lem}\label{lem12}
The group $G$ is a subgroup of $H$.
\end{lem}

{\bf Proof.}
Take $\lambda\in {\bf C}$ and a projection $p$ from
$M_{\infty}(A)$. We have to demonstrate that
$g_{\lambda ;p}^{(k)}$ belongs to $H$ for all $k\ge 1$.
Let us prove this statement by induction over $k$.
The case $k=1$ is clear. Suppose
$g_{\lambda ;p}^{(k)}\in H$ for all $k\le n-1$.
In particular,
$$
g_{\lambda ;p}^{(n-1)}=
[\lambda p^{\oplus (n-1)}]-[(n-1)\lambda p]=y
$$
for some $y\in H$. Therefore,
$$
\begin{array}{ccl}
g_{\lambda ;p}^{(n)} &=&
[\lambda p^{\oplus (n-1)}]+[\lambda p]-
[((n-1)\lambda+\lambda)p]\\
&=& [(n-1)\lambda p]+y+[\lambda p]-[((n-1)\lambda+\lambda)p]\\
&=& y-h_{(n-1)\lambda,\lambda;p}.
\end{array}
$$
Thus $g_{\lambda ;p}^{(n)}\in H$ and by induction
we obtain the desired statement.
$\pcom$

Let us consider the map
\begin{equation}\label{eq10}
h:N_0(A)_{fin}\lra K_0(A)\otimes{\bf C}
\end{equation}
defined as follows
$$ h([a]-[b])=
\sum_{i=1}^n [P_a(\lambda_i)]\otimes\lambda_i-
\sum_{j=1}^m[P_b(\mu_j)]\otimes\mu_j,$$
where $a,b\in M_{\infty}(A)_{fin}$ and
$sp(a)=\{\lambda_1,\dots,\lambda_n\}$, $sp(b)=\{\mu_1,\dots,\mu_m\}$.

\begin{prop}
The map $h$ is a surjective homomorphism of groups.
Besides, the kernel of $h$ coinsides with the group $H$.
\end{prop}

{\bf Proof.}
It is obvious that $h$ is a well defined surjective homomorphism.
Also, it is clear that $H$ belongs to the kernel of $h$.
To complete the proof let us construct the inverse for $h$
homomorphism
$$
t:K_0(A)\otimes{\bf C}\lra N_0(A)_{fin}/H.
$$
We put
$$
t(\sum_{i=1}^n ([p_i]-[q_i])\otimes\lambda_i)=
\sum_{i=1}^n[p_i\lambda_i]-\sum_{i=1}^n[q_i\lambda_i]+H.
$$
Let us demonstrate that the map $t$ is well defined.
In the other words, we have to verify that the homomorphism
$t$ is trivial on the elements
$c_{\lambda,\mu ;p}^{(1)}=
[p]\otimes (\lambda +\mu)-[p]\otimes\lambda -[p]\otimes\mu$,
$c_{z;\lambda ;p}^{(2)}=
[p]z\otimes\lambda -[p]\otimes z\lambda$,
and
$c_{\lambda ;p,q}^{(3)}=
([p]+[q])\otimes\lambda -[p]\otimes\lambda -[q]\otimes\lambda$,
where $\lambda ,\mu\in{\bf C}, z\in{\bf Z}$ and
$p,q\in M_{\infty}(A)$ are projections.
We derive
$t(c_{\lambda,\mu ;p}^{(1)})=h_{\lambda ,\mu ;p}+H=H$.
Besides, $c_{-1;\lambda ;p}^{(2)}=
c_{\lambda,-\lambda ;p}^{(1)}$.
Therefore it suffices to regard the elements
$c_{z;\lambda ;p}^{(2)}$ provided $z>0$.
In this case we conclude
$t(c_{z;\lambda ;p}^{(2)})=
t([p^{\oplus z}]\otimes\lambda-[p]\otimes z\lambda)=
g_{\lambda ;p}^{(z)}+H=H$,
where we have used Lemma \ref{lem12}.
Finally, it can be directly verified that
$t(c_{\lambda ;p,q}^{(3)})=0+H=H$.
To complete the proof, it remains to note that $t=h^{-1}$.
$\pcom$

\section{Banach cyclic homology}
As above, let $A$ be a von Neumann algebra.
First let us recall some concepts from noncommutative geometry
(see, for example,~\cite{Con,Lod,Sol-Tro}).
Consider the complex vector space $C_n(A)=A^{\otimes (n+1)}$, where
$A^{\otimes (n+1)}=\underbrace{A\otimes A\otimes \dots\otimes A}_{n+1}$.
The {\em cyclic operator\/} $\tau_n: C_n(A)\lra C_n(A)$
is defined on generators by the formula
$$\tau_n(a_0\otimes\dots\otimes a_n)={(-1)}^n a_n\otimes a_0
\otimes\dots\otimes a_{n-1}.$$
The cokernel of the endomorphism $1-\tau_n: C_n(A)\lra C_n(A)$
we denote by
$$
CC_n(A)=A^{\otimes (n+1)}/\im (1-\tau_n).
$$
Further, we define the {\em face operator\/}
$b_n:C_n(A)\lra C_{n-1}(A)$ by the formula
$$
b_n(a_0\otimes a_1\otimes\dots\otimes a_n)=
  \sum_{i=0}^{n-1}{(-1)}^i a_0\otimes a_1\otimes\dots\otimes
   a_ia_{i+1}\otimes\dots\otimes a_n+ $$
$$
{(-1)}^n a_na_0\otimes a_1\otimes\dots\otimes a_{n-1}.
$$
It is clear that
\begin{equation}\label{eq1}
b=\sum_{i=0}^n {(-1)}^i d_i,
\end{equation}
where the linear maps
$d_i:A^{\otimes (n+1)}\lra A^{\otimes n}$
are defined as follows
$$
\begin{array}{lcc}
d_i (a_0\otimes\dots\otimes a_n) &=&
a_0\otimes\dots\otimes a_ia_{i+1}\otimes\dots\otimes a_n,\\
&& 0\le i\le n-1, \\
d_n (a_0\otimes\dots\otimes a_n) &=&
a_na_0\otimes a_1\otimes\dots\otimes a_{n-1}.
 \end{array}
$$

It can be verified by the direct calculation
that the family of linear spaces
$CC_*(A)=\{CC_n(A),b_n\}$ is a chain complex.
Homology of this complex is called {\em cyclic homology\/} of $A$
and is denoted by $HC_n(A)=H_n(CC_*(A))$, $n\ge 0$.

The {\em trace map\/}
$Tr: M_r(A)^{\otimes (n+1)}\lra A^{\otimes (n+1)}$
is defined by the formula
$$
Tr(\xi^{(0)}\otimes\xi^{(1)}\otimes\dots\otimes\xi^{(n)})=
\sum_{i_0,\dots,i_n=1}^r
\xi^{(0)}_{i_0i_1}\otimes\xi^{(1)}_{i_1i_2}\otimes\dots\otimes
\xi^{(n)}_{i_ni_0},
$$
where $\xi^{(k)}=(\xi^{(k)}_{i,j})_{i,j=1}^r\in M_r(A)$.
It can be directly verified that the trace map
$Tr: CC_*(M_r(A))\lra CC_*(A)$
is a morphism of chain complexes.
Furthermore, the induced map
$$
Tr_*: HC_*(M_r(A))\stackrel{\cong}{\lra} HC_*(A)
$$
is an isomorphism.

Now let $X,Y$ be normed spaces, and $x\in X, y\in Y$.
A representative of the equivalence class
$x\otimes y\in X\otimes Y$ we shall denote by
$x\Box y$. Assume
$$
\|\sum_i x_i\Box y_i\|:=\sum_i \|x_i\|\|y_i\|
\quad (x_i\in X, y_i\in Y).
$$
Then the {\em projective norm\/}
of an equivalence class $\xi\in X\otimes Y$ is defined as follows
$$
\|\xi\|=\inf \{\|\sum_i x_i\Box y_i\| :
x_i\in X, y_i\in Y \;\mbox{and}\; \sum_i x_i\Box y_i \in\xi \}.
$$

Below we assume that all tensor products of normed spaces
are equipped with the projective norm.

Under the previous conventions, let us set
$$
{\cal CC}_n(A)=A^{\otimes (n+1)}/\overline{\im (1-\tau_n)}.
$$
Note that ${\cal CC}_n(A)$ is a Banach space.
For $\xi\in A^{\otimes (n+1)}$ we denote by
$[\xi]_{CC_n(A)}$ the quotient class of $\xi$ in $CC_n(A)$,
and by $[\xi]_{{\cal CC}_n(A)}$ the quotient class of $\xi$ in
${\cal CC}_n(A)$.

\begin{lem}\label{lem3}
The face operator $b_n:A^{\otimes (n+1)}\lra A^{\otimes n}$
is a continuous map.
\end{lem}

{\bf Proof.}
Under equality (\ref{eq1}), it is sufficiently to prove
that all maps $d_i$ $(0\le i\le n)$ are continuous.
Given any $\varepsilon >0$. For each $\xi\in A^{\otimes (n+1)}$
we can find an element
$\sum_k a_0^{(k)}\Box\dots\Box a_n^{(k)}\in\xi$
such that
$$\|\sum_k a_0^{(k)}\Box\dots\Box a_n^{(k)}\|-\|\xi\|
<\varepsilon.$$
Therefore,
$$
\begin{array}{lcl}
\|d_i (\xi)\| &=&
\|\sum_k d_i(a_0^{(k)}\otimes\dots\otimes a_n^{(k)})\|\\
&=& \|\sum_k a_0^{(k)}\otimes\dots\otimes
a_i^{(k)}a_{i+1}^{(k)}\otimes\dots\otimes a_n^{(k)}\| \\
&\le& \|\sum_k a_0^{(k)}\Box\dots\Box
a_i^{(k)}a_{i+1}^{(k)}\Box\dots\Box a_n^{(k)}\|  \\
&=& \sum_k \|a_0^{(k)}\|\dots
\|a_i^{(k)}a_{i+1}^{(k)}\|\dots \|a_n^{(k)}\|  \\
&\le& \sum_k \|a_0^{(k)}\|\dots
\|a_i^{(k)}\|\|a_{i+1}^{(k)}\|\dots \|a
_n^{(k)}\|\\
&=& \|\sum_k a_0^{(k)}\Box\dots\Box a_n^{(k)}\|
<\|\xi\|+\varepsilon
\end{array}
$$
for all $\varepsilon>0$. Hence, $\|d_i (\xi)\|\le\|\xi\|$.
$\pcom$

Let us define the map
$\beta_n:{\cal CC}_n(A)\lra {\cal CC}_{n-1}(A)$
by the formula
$$
\beta_n ([x]_{{\cal CC}_n(A)})=[b_n (x)]_{{\cal CC}_{n-1}(A)}.
$$
From Lemma~\ref{lem3} we conclude that the family
${\cal CC}_*(A)=\{{\cal CC}_n(A),\beta_n\}$ is a
well defined chain complex.
Let us put
$$
{\cal HC}_n(A)=
\Ker\beta_n/ \,\overline{ \im\beta_{n+1}}.
$$
The quotient space ${\cal HC}_*(A)$ we shall call
{\em Banach (cyclic) homology\/} of $A$
(cf.~\cite{Hel1,Hel2}).
Note that  ${\cal HC}_n(A)$ is a Banach space (for each $n\ge 0$).

For $\xi\in M_r(A)^{\otimes (n+1)}$ let us denote by
$\langle\xi\rangle_{HC}=\langle\xi\rangle_{HC_*(M_r(A))}$
the cyclic homology class of $\xi$ and by
$\langle\xi\rangle_{{\cal HC}}=
\langle\xi\rangle_{{\cal HC}_*(M_r(A))}$
the Banach cyclic homology class of $\xi$.

Let $p\in M_r(A)$ be projection. Then
$$
b_{2l}(p^{\otimes (2l+1)})=\sum_{k=0}^{2l}
d_k (p^{\otimes (2l+1)})=
\sum_{k=0}^{2l} {(-1)}^k p^{\otimes 2l}
=p^{\otimes 2l}.
$$

On the other hand, we have
$[p^{\otimes 2l}]_{CC_*(M_r(A))}
={(-1)}^{(2l-1)}[p^{\otimes 2l}]_{CC_*(M_r(A))}$.
So $[p^{\otimes 2l}]_{CC_*(M_r(A))}=0$.
Therefore $p^{\otimes (2l+1)}$ is a cycle.
Now let us define the {\em Chern character\/}
$$ Ch^0_{2l}:K_0(A)\lra HC_{2l}(A)$$ by
$Ch^0_{2l}([p])=Tr_*\langle {(-1)}^l p^{\otimes (2l+1)}\rangle_{HC}$.
The Chern character is a well defined homomorphism of
groups~\cite[Theorem 8.3.2]{Lod}.

Let us study the linear epimorphism
$\pi_n :CC_n(A)\lra {\cal CC}_n(A)$, where
$$
\pi_n([\xi]_{CC_n(A)})=[\xi]_{{\cal CC}_n(A)},
\;\xi\in A^{\otimes (n+1)}.
$$

It is obvious that the family of the maps
$\{\pi_n\}: CC_*(A)\lra {\cal CC}_*(A)$
is a chain homomorphism, i.e.,
$\pi_{n-1} b_n=\beta_n\pi_n$ for all $n\ge 1$.
So $\pi_n (\Ker b_n)\subset \Ker \beta_n$ and
$\pi_n (\im b_{n+1})\subset \overline{\im \beta_{n+1}}$.
Therefore,
\begin{equation}\label{eq2}
\pi_*:HC_*(A)\lra {\cal HC}_*(A)
\end{equation}
is a well defined map of homology spaces.

Below we shall need the following result.
\begin{lem}\label{lem5}
The trace $Tr: M_r(A)^{\otimes (n+1)}\lra A^{\otimes (n+1)}$
is a continuous map.
\end{lem}

{\bf Proof.}
Given any $\varepsilon >0$. For each quotient class
$\xi\in {M_r(A)}^{\otimes (n+1)}$ of the tensor product
we can find a representative
$\sum_{k=1}^N
\xi^{(0),k}\Box\xi^{(1),k}\Box\dots\Box\xi^{(n),k}$ of $\xi$
such that
$$\|\sum_{k=1}^N
\xi^{(0),k}\Box\xi^{(1),k}\Box\dots\Box\xi^{(n),k}\|-
\|\xi\|<\varepsilon.$$

Therefore,
$$
\begin{array}{lcl}
\|Tr (\xi)\| &=&
\|\sum_{k=1}^N Tr(\xi^{(0),k}\otimes\xi^{(1),k}
\otimes\dots\otimes\xi^{(n),k})\| \\
&=& \|\sum_{k=1}^N
\sum_{i_0,\dots,i_n=1}^r
\xi^{(0),k}_{i_0i_1}\otimes\xi^{(1),k}_{i_1i_2}\otimes\dots\otimes
\xi^{(n),k}_{i_ni_0}\| \\
&\le& \|\sum_{k=1}^N \sum_{i_0,\dots,i_n=1}^r
\xi^{(0),k}_{i_0i_1}\Box\xi^{(1),k}_{i_1i_2}\Box\dots\Box
\xi^{(n),k}_{i_ni_0}\| \\
&=& \sum_{k=1}^N \sum_{i_0,\dots,i_n=1}^r
\|\xi^{(0),k}_{i_0i_1}\|\|\xi^{(1),k}_{i_1i_2}\|\dots
\|\xi^{(n),k}_{i_ni_0}\| \\
&\le& \sum_{k=1}^N \sum_{i_0,\dots,i_n=1}^r
\|\xi^{(0),k}\|\|\xi^{(1),k}\|\dots\|\xi^{(n),k}\| \\
&=& r^{n+1} \|\sum_{k=1}^N
\xi^{(0),k}\Box\xi^{(1),k}\Box\dots\Box\xi^{(n),k}\|
< r^{n+1}(\|\xi\|+\varepsilon )
\end{array}
$$
for all $\varepsilon>0$. So, $\|Tr (\xi)\|\le r^{n+1}\|\xi\|$.
$\pcom$

From Lemma~\ref{lem5} we conclude that the map
\begin{equation}\label{eq9}
Tr :{\cal CC}_*(M_r(A))\lra {\cal CC}_*(A),\;
[\xi]_{{\cal CC}_*(M_r(A))}\longmapsto [Tr(\xi)]_{{\cal CC}_*(A)}
\end{equation}
is well defined.
Furthermore, it can be directly checked up that
trace map~(\ref{eq9}) is a morphism of chain complexes.
By Lemma~\ref{lem5}, this implies that the induced homomorphism
$$Tr_*:{\cal HC}_*(M_r(A))\lra {\cal HC}_*(A)$$
of Banach homology is well defined.

\section{The generalized Chern character}
Given $l\ge 0$. We want to define a map
$$
T:M_{\infty}(A)_{\nu}\lra {\cal HC}_{2l}(A)
$$
to the even Banach homology by the following construction.
Let an element $a$ belong to $M_\infty(A)_{\nu}$.
Then we can suppose that $a\in M_r(A)$ for some $r\ge 1$.
For each natural number $n$ let us consider a cover
${\cal E}^{(n)}=\{E_k^{(n)}\}_{k=1}^{k_n}$
of the spectrum of $a$ by disjoint Borel sets such that
the diameter of each of these sets does not exceed $1/n$.
Moreover, we can assume that
$\{E_k^{(m)}\}_{k=1}^{k_m}$ is a subdivision
of the cover $\{E_k^{(n)}\}_{k=1}^{k_n}$ when $m\ge n$.
Therefore we can write
$$
E_k^{(n)}=\bigsqcup_{j=1}^{j(k)} E_{k,j}^{(m)},
$$
where $E_{k,j}^{(m)}$ are some elements of the cover ${\cal E}^{(m)}$.
Thus, ${\cal E}^{(m)}=\{E_{k,j}^{(m)}\}_{k=1;j=1}^{k_n;j(k)}$ and
$\sum_{k=1}^{k_n}j(k)=k_m$.

Also, for any $\lambda_k^{(n)}\in E_k^{(n)}$ let us consider
a sequence $\{a_n\}_{n=1}^\infty$ from $M_r(A)$, where
$$
a_n=\sum_{k=1}^{k_n} P_a(E_k^{(n)})\lambda_k^{(n)}.
$$
It follows from the spectral theorem that
$\{a_n\}_{n=1}^\infty$ converges
uniformly to $a$.
Further, for each $a_n$ let us examine
$$
\tilde a_n=\sum_{k=1}^{k_n}P_a(E_k^{(n)})^{\otimes (2l+1)}
 \lambda_k^{(n)}\in M_r(A)^{\otimes (2l+1)}.
$$

Given natural numbers $n,m$ $(m\ge n)$.
Under the previous notation, we see that
$$
a_n=\sum_{k=1}^{k_n}P_a(\bigsqcup_{j=1}^{j(k)} E_{k,j}^{(m)})
\lambda_k^{(n)}
=\sum_{k=1}^{k_n}\sum_{j=1}^{j(k)}
P_a(E_{k,j}^{(m)})\lambda_k^{(n)}
$$
and $a_m=\sum_{k=1}^{k_n}\sum_{j=1}^{j(k)} P_a(E_{k,j}^{(m)})
\lambda_{k,j}^{(m)}$,
where $\lambda_{k,j}^{(m)}\in E_{k,j}^{(m)}$.
This yields that
\begin{equation}\label{eq3}
a_n-a_m=\sum_{k=1}^{k_n}\sum_{j=1}^{j(k)} P_a(E_{k,j}^{(m)})
(\lambda_k^{(n)}-\lambda_{k,j}^{(m)}).
\end{equation}

Furthermore,
\begin{eqnarray*}
\tilde a_n&=&\sum_{k=1}^{k_n}P_a(E_k^{(n)})^{\otimes (2l+1)}
\lambda_k^{(n)}
=\sum_{k=1}^{k_n} P_a(\bigsqcup_{j=1}^{j(k)}
E_{k,j}^{(m)})^{\otimes (2l+1)} \lambda_k^{(n)} \\
&=&\sum_{k=1}^{k_n}(\sum_{j_0=1}^{j(k)}P_a(E_{k,j_0}^{(m)}))
\otimes\dots\otimes (\sum_{j_{2l}=1}^{j(k)}P_a(E_{k,j_{2l}}^{(m)}))
\lambda_k^{(n)}  \\
&=&\sum_{k=1}^{k_n}(\sum_{j_0=1}^{j(k)}\dots
\sum_{j_{2l}=1}^{j(k)}
P_a(E_{k,j_0}^{(m)}) \otimes\dots\otimes P_a(E_{k,j_{2l}}^{(m)}))
\lambda_k^{(n)}
\end{eqnarray*}
and
$\tilde a_m=\sum_{k=1}^{k_n}\sum_{j=1}^{j(k)}
P_a(E_{k,j}^{(m)})^{\otimes (2l+1)}\lambda_{k,j}^{(m)}$.

For brevity we shall use the following notation
$${\widetilde{\sum}}_{j_0,\dots,j_{2l}=1}^{j(k)} :=
\sum_{\{1\le j_0\le\dots\le j_{2l}\le j(k) \,:\,\exists j_p\ne j_q\}}
$$
for the sum over $j_0,\dots,j_{2l}$ from one to $j(k)$,
where not all indices coinside.

So we obtain
\begin{eqnarray}\label{eq4}
\tilde a_n-\tilde a_m &=& \sum_{k=1}^{k_n}\sum_{j=1}^{j(k)}
P_a(E_{k,j}^{(m)})^{\otimes (2l+1)}
(\lambda_k^{(n)}-\lambda_{k,j}^{(m)}) \nonumber\\
&+& \sum_{k=1}^{k_n}
({\widetilde{\sum}}^{j(k)}_{j_0,\dots,j_{2l}=1}
P_a(E_{k,j_0}^{(m)}) \otimes\dots\otimes P_a(E_{k,j_{2l}}^{(m)}))
\lambda_k^{(n)}\\
&=& \alpha_{n,m}+\gamma_{n,m}, \nonumber
\end{eqnarray}
where by $\alpha_{n,m}$ $(\gamma_{n,m})$ we denote the first
(the second) summand in expression~(\ref{eq4}).

The following result is the main ingredient of our definition
of the map~$T$.

\begin{th}\label{th2}
Let $\{p_i\}_{i=1}^N\in M_r(A)$ be a family of pairwise
orthogonal projections and
$\eta ={\widetilde{\sum}}_{j_0,\dots,j_{2l}=1}^N
p_{j_0}\otimes\dots\otimes p_{j_{2l}}$.
Then the element $\eta$ belongs to the kernel of the face
operator $\beta_{2l}$. Besides,
$\langle \eta \rangle_{{\cal HC}_{2l}(M_r(A))} =0$.
\end{th}

{\bf Proof.}
We have
$$
\eta=(\sum_{j=1}^N p_j)^{\otimes (2l+1)}-
\sum_{j=1}^N p_j^{\otimes (2l+1)}.
$$
So $\eta$ belongs to the kernel of the face operator $\beta_{2l}$
as a difference of elements from $\Ker\beta_{2l}$.

Now let us examine the following element
$$
\begin{array}{crl}
\alpha &:=& Tr_* \langle
(\sum_{j=1}^N p_j)^{\otimes (2l+1)}\rangle_{HC}\\
&=& {(-1)}^l Ch^0_{2l} (\sum_{j=1}^N [p_j])=
\sum_{j=1}^N {(-1)}^l Ch^0_{2l} ([p_j])\\
&=& \sum_{j=1}^N Tr_* \langle
p_j^{\otimes (2l+1)}\rangle_{HC}
=Tr_* \langle \sum_{j=1}^N p_j^{\otimes (2l+1)}\rangle_{HC}.
\end{array}
$$

On the other hand,
$$
\alpha =Tr_*\langle \sum_{j_0,\dots,j_{2l}=1}^N
p_{j_0}\otimes\dots\otimes p_{j_{2l}}\rangle_{HC}.
$$
Thus,
$$
\begin{array}{ccl}
0 &=& Tr_*\langle \sum_{j_0,\dots,j_{2l}=1}^N
p_{j_0}\otimes\dots\otimes p_{j_{2l}}\rangle_{HC}
-Tr_* \langle \sum_{j=1}^N p_j^{\otimes (2l+1)}\rangle_{HC}\\
&=&Tr_*\langle \sum_{j_0,\dots,j_{2l}=1}^N
p_{j_0}\otimes\dots\otimes p_{j_{2l}}
-\sum_{j=1}^N p_j^{\otimes (2l+1)}\rangle_{HC}\\
&=& Tr_*\langle
{\widetilde{\sum}}_{j_0,\dots,j_{2l}=1}^N
p_{j_0}\otimes\dots\otimes p_{j_{2l}}\rangle_{HC}.
\end{array}
$$

The trace map is an isomorphism. Therefore we conclude
$$\langle\eta\rangle_{HC}
=\langle {\widetilde{\sum}}_{j_0,\dots,j_{2l}=1}^N
p_{j_0}\otimes\dots\otimes p_{j_{2l}}\rangle_{HC}=0. $$
Whence, $\langle \eta \rangle_{{\cal HC}}=
\pi_*(\langle\eta\rangle_{HC})=0$,
where $\pi_*$ is map~(\ref{eq2}).
$\pcom$

Now let us return to expressions~(\ref{eq3}),(\ref{eq4}).
By the previous theorem, we see that
$\langle\gamma_{n,m}\rangle_{{\cal HC}}=0$
so $\langle\tilde a_n-\tilde a_m\rangle_{{\cal HC}}=
\langle\alpha_{n,m}\rangle_{{\cal HC}}$.

We claim that $\|\alpha_{n,m}\|=\|a_n-a_m\|$. Indeed,
it is clear that elements $a_n-a_m$ and $\alpha_{n,m}$
are normal. Besides, $sp (a_n-a_m)\setminus\{0\}$ and
$sp (\alpha_{n,m})\setminus\{0\}$
coinside with the set
$\{\lambda_k^{(n)}-\lambda_{k,j}^{(m)}
\}_{k=1; j=1}^{k_n; j(k)}$.
This implies that spectral radii of these elements coinside too.
Thus we obtain the desired statement.

So we can write
\begin{equation}\label{eq5}
\|\langle\tilde a_n-\tilde a_m\rangle_{{\cal HC}}\|
=\|\langle\alpha_{n,m}\rangle_{{\cal HC}}\|
\le \|\alpha_{n,m}\|=\|a_n-a_m\|.
\end{equation}

By Lemma~\ref{lem5} and inequality~(\ref{eq5}),
we conclude that
$\{Tr_*(\langle\tilde a_n\rangle_{{\cal HC}})\}_{n=1}^\infty$
is a Cauchy sequence. Therefore it converges to some element
$T(a;\{a_n\})\in {\cal HC}_{2l}(A)$.

It remains to verify that the limit $T(a;\{a_n\})$
does not depend on $\{a_n\}$. Let us regard covers
${\cal E}^{(n)}=\{E_k^{(n)}\}_{k=1}^{k_n}$ and
${\cal F}^{(n)}=\{F_j^{(n)}\}_{j=1}^{j_n}$
of the spectrum of $a$ by disjoint Borel sets.
Besides, we shall assume that the diameter of each
of these sets does not exceed $1/n$.
Also, for any $\mu_j^{(n)}\in F_j^{(n)}$
let us examine a sequence $\{c_n\}_{n=1}^\infty$,
where $c_n=\sum_{j=1}^{j_n} P_a(F_j^{(n)})\mu_j^{(n)}$.

\begin{th}\label{th6}
The elements $T(a;\{a_n\})$ and $T(a;\{c_n\})$ coinside.
Thus the map
$$T:M_{\infty}(A)_{\nu}\lra {\cal HC}_{2l}(A),\quad
a\longmapsto T(a;\{a_n\})=T(a)$$
is well defined.
\end{th}

{\bf Proof.}
Let us denote $X_{k,j}^{(n)}=E_k^{(n)}\cap F_j^{(n)}$.
Then
$$
 \begin{array}{rcl}
a_n =\sum_{k=1}^{k_n} P_a(E_k^{(n)})\lambda_k^{(n)}
&=&\sum_{k=1}^{k_n}P_a(\bigsqcup_{j=1}^{j(n)} E_k^{(n)}\cap
F_j^{(n)})\lambda_k^{(n)}\\
&=&\sum_{k=1}^{k_n}\sum_{j=1}^{j_n} P_a(X_{k,j}^{(n)})
\lambda_k^{(n)}
 \end{array}
$$
and by the same reason
$c_n=\sum_{k=1}^{k_n}\sum_{j=1}^{j_n} P_a(X_{k,j}^{(n)})
\mu_j^{(n)}$.
Hence,
\begin{equation}\label{eq6}
a_n-c_n=\sum_{k=1}^{k_n}\sum_{j=1}^{j_n}P_a(X_{k,j}^{(n)})
(\lambda_k^{(n)}-\mu_j^{(n)}).
\end{equation}
If $X_{k,j}^{(n)}=\emptyset$, then $P_a(X_{k,j}^{(n)})=0$.
Therefore we can assume that
$X_{k,j}^{(n)}\ne\emptyset$ in expression (\ref{eq6}).
In this case let us consider $z\in X_{k,j}^{(n)}$.
We deduce that
$|\lambda_k^{(n)}-\mu_j^{(n)}|\le
|\lambda_k^{(n)}-z|+|z-\mu_j^{(n)}|\le 1/n+1/n=2/n$
for all $1\le k\le k_n$, $1\le j\le j_n$.
Note that $a_n-c_n$ is a normal element.
Therefore,
\begin{equation}\label{eq7}
\|a_n-c_n\|= \sup\{|\lambda_k^{(n)}-\mu_j^{(n)}| :
X_{k,j}^{(n)}\ne\emptyset, 1\le k\le k_n, 1\le j\le j_n\}
\le 2/n.
\end{equation}

On the other hand, we have
$$
\begin{array}{ccl}
\tilde a_n &=&\sum_{k=1}^{k_n}P_a(E_k^{(n)})^{\otimes (2l+1)}
\lambda_k^{(n)}
=\sum_{k=1}^{k_n}(P_a(\bigsqcup_{j=1}^{j(n)}
X_{k,j}^{(n)}))^{\otimes (2l+1)} \lambda_k^{(n)}\\
&=& \sum_{k=1}^{k_n}\sum_{j_0,\dots j_{2l}=1}^{j_n}
P_a(X_{k,j_0}^{(n)})\otimes\dots\otimes P_a(X_{k,j_{2l}}^{(n)})
\lambda_k^{(n)}
\end{array}
$$
and by the same reason
$$
\tilde c_n=\sum_{j=1}^{j_n}\sum_{k_0,\dots k_{2l}=1}^{k_n}
P_a(X_{j,k_0}^{(n)})\otimes\dots\otimes P_a(X_{j,k_{2l}}^{(n)})
\mu_j^{(n)}.
$$
Thus we obtain
$$
\begin{array}{ccl}
\tilde a_n-\tilde c_n &=& \sum_{k=1}^{k_n}
\sum_{j=1}^{j_n} P_a(X_{k,j}^{(n)})^{\otimes (2l+1)}
(\lambda_k^{(n)}-\mu_j^{(n)}) \\
&+& \sum_{k=1}^{k_n}
({\widetilde{\sum}}_{j_0,\dots,j_{2l}=1}^{j(n)}
P_a(X_{k,j_0}^{(n)})\otimes\dots\otimes P_a(X_{k,j_{2l}}^{(n)}))
\lambda_k^{(n)}\\
&-& \sum_{j=1}^{j_n}
({\widetilde{\sum}}_{k_0,\dots,k_{2l}=1}^{k(n)}
P_a(X_{j,k_0}^{(n)})\otimes\dots\otimes P_a(X_{j,k_{2l}}^{(n)}))
\mu_j^{(n)} \\
&=& \gamma_n^{(1)}+\gamma_n^{(2)}-\gamma_n^{(3)}.
\end{array}
$$

From Theorem~\ref{th2} we conclude that
$\langle\gamma_n^{(2)}\rangle_{{\cal HC}}
=\langle\gamma_n^{(3)}\rangle_{{\cal HC}}=0$.
The elements $a_n-c_n$ and $\gamma_n^{(1)}$ are normal.
Furtermore, the sets $sp(a_n-c_n)\setminus\{0\}$ and
$sp(\gamma_n^{(1)})\setminus\{0\}$ coinside.
Therefore, $\|a_n-c_n\|=\|\gamma_n^{(1)}\|$.
Using inequality~(\ref{eq7}), we obtain
$$\|\langle\tilde a_n\rangle_{{\cal HC}}-
\langle\tilde c_n\rangle_{{\cal HC}}\|
=\|\langle\gamma_n^{(1)}\rangle_{{\cal HC}}\|
\le\|\gamma_n^{(1)}\|=\|a_n-c_n\|\le 2/n.$$
Therefore,
$\|Tr_*(\langle\tilde a_n\rangle_{{\cal HC}})-
Tr_*(\langle\tilde c_n\rangle_{{\cal HC}})\|\le
2\|Tr_*\|/n$.
This estimate implies that
$$
T(a,\{a_n\})
=\lim_n Tr_*(\langle\tilde a_n\rangle_{{\cal HC}})
=\lim_n Tr_*(\langle\tilde c_n\rangle_{{\cal HC}})
=T(a,\{c_n\}).
$$
The proof is complete.
$\pcom$

\begin{prop}\label{lem11}
Suppose $a,b\in M_\infty(A)_{\nu}$ are equivalent in the sence of
Definition~{\rm\ref{dfn1}}. Then $T(a)=T(b)$.
\end{prop}

{\bf Proof.}
For each $n\in {\bf N}$ let us cover the space $sp(a)\cap sp(b)$
by disjoint Borel sets and enlarge this system of sets
to a disjoint cover
$\{E_k^{(n)}\}_{k=0}^{k_n}$ of the space $sp(a)\cup sp(b)$.
As above, we suppose that
${\rm diam}(E_k^{(n)})\le 1/n$ for all $0\le k\le k_n$.
Also, let us consider $\lambda_k^{(n)}\in E_k^{(n)}$.
If $0\in sp(a)\cup sp(b)$, then we shall assume that $0\in E_0^{(n)}$
and $\lambda_0^{(n)}=0$. In the opposite case,
we put $E_0^{(n)}=\emptyset$.
Furtemore, for any $1\le k\le k_n$ we can assume that
$E_k^{(n)}$ is an admissible Borel set.

Projections $P_a(E_k^{(n)})$ and $P_b(E_k^{(n)})$
are stably equivalent for all ${k,n\ge 1}$. Therefore,
$Ch_{2l}^0 ([P_a(E_k^{(n)})])=Ch_{2l}^0 ([P_b(E_k^{(n)})])$.
Since the trace map is an isomorphism, we conclude that
$\langle {P_a(E_k^{(n)})}^{\otimes (2l+1)}\rangle_{HC}=
\langle {P_b(E_k^{(n)})}^{\otimes (2l+1)}\rangle_{HC}$.
Thus,
$$
\begin{array}{ccl}
\langle {P_a(E_k^{(n)})}^{\otimes (2l+1)}\rangle_{{\cal HC}}
&=&\pi_* (\langle {P_a(E_k^{(n)})}^{\otimes (2l+1)}\rangle_{HC})\\
&=&\pi_* (\langle {P_b(E_k^{(n)})}^{\otimes (2l+1)}\rangle_{HC})
=\langle {P_b(E_k^{(n)})}^{\otimes (2l+1)}\rangle_{{\cal HC}},
\end{array}
$$
where $\pi_*$ is map~(\ref{eq2}). So we obtain that
$$
\langle\tilde a_n\rangle_{{\cal HC}}
=\sum_{k=0}^{k_n} \langle {P_a(E_k^{(n)})}^{\otimes (2l+1)}
\rangle_{{\cal HC}}\lambda_k^{(n)}
=\sum_{k=0}^{k_n} \langle {P_b(E_k^{(n)})}^{\otimes (2l+1)}
\rangle_{{\cal HC}}\lambda_k^{(n)}
=\langle\tilde b_n\rangle_{{\cal HC}}
$$
for all $n\ge 1$ so $T(a)=T(b)$.
$\pcom$

\begin{dfn}
{\rm We define the}  generalized Chern character
{\rm as the map }
$$
{\cal C}h_{2l}^0 : N_0(A)\lra {\cal HC}_{2l}(A),\quad
[a]-[b]\longmapsto {(-1)}^l (T(a)-T(b)).
$$
\end{dfn}

It follows from Proposition~\ref{lem11} that the generalized
Chern character is well defined. An immediate verification gives us

 \begin{prop}
The generalized Chern character is a homomorphism of
groups.~$\pcom$
 \end{prop}

\begin{th}\label{th7}
For any $l\ge 0$ there is a commutative diagramm
$$
\begin{array}{ccc}
K_0(A) &\hookrightarrow &N_0(A) \\
\downarrow \lefteqn{Ch_0^{2l}} && \downarrow
\lefteqn{{\cal C}h_0^{2l}}\\
HC_{2l}(A) &\stackrel{\pi_*}{\lra}& {\cal HC}_{2l}(A),
\end{array}
$$
where $\pi_*$ is map~{\rm (\ref{eq2})}.
\end{th}

{\bf Proof.}
Under the notation of the beginning of this section
let us argue as follows.
Let $p\in M_\infty (A)$ be a projection.
In this case the cover ${\cal E} ^{(n)}$ of the spectrum of $p$
coinside with the set $\{\{1\},\{0\}\}$ for all $n\ge 1$.
Therefore for all $n\ge 1$ we have $\tilde p_n =p^{\otimes (2l+1)}$.
Finally, we obtain
$$
\begin{array}{ccl}
{\cal C}h_0^{2l}([p]) &=&
\lim_n Tr_*({(-1)}^l\langle\tilde p_n\rangle_{{\cal HC}_{2l}})\\
&=& \lim_n Tr_*({(-1)}^l\langle p^{\otimes (2l+1)}
\rangle_{{\cal HC}_{2l}})\\
&=& Tr_*({(-1)}^l\langle p^{\otimes (2l+1)}
\rangle_{{\cal HC}_{2l}})\\
&=& \langle{(-1)}^l Tr(p^{\otimes (2l+1)})
\rangle_{{\cal HC}_{2l}}\\
&=& \pi_* (\langle{(-1)}^l Tr(p^{\otimes (2l+1)})
\rangle_{HC_{2l}})
= \pi_* Ch_{2l}^0([p]).
\end{array}
$$
The proof is complete.
$\pcom$

\begin{th}\label{th8}
For any $l\ge 0$ there is a commutative diagramm
$$
\begin{array}{ccc}
N_0(A)_{fin} &\stackrel{h}{\lra} &K_0(A)\otimes {\bf C} \\
\downarrow \lefteqn{{\cal C}h_0^{2l}} && \downarrow
\lefteqn{\widetilde{Ch}_0^{2l}}\\
{\cal HC}_{2l}(A) &\stackrel{\pi_*}{\longleftarrow}&
HC_{2l}(A),
\end{array}
$$
where $h$ is map~{\rm (\ref{eq10})} and
$\widetilde{Ch}_{2l}^0([p]\otimes\lambda)=
Ch_{2l}^0([p])\lambda$.
\end{th}

{\bf Proof.}
Let $a$ be an element of $M_\infty(A)_{fin}$ such that
$sp(a)=\{\lambda_1,\dots ,\lambda_n\}$.
Assume $p_i:=P_a(\{\lambda_i\})$. Then we have
$$
\begin{array}{ccl}
\pi_* \widetilde{Ch}_0^{2l} h ([a]) &=&
\pi_* \widetilde{Ch}_0^{2l}
(\sum_{i=1}^n [p_i]\otimes\lambda_i)\\
&=&\pi_* (\sum_{i=1}^n Ch_0^{2l} ([p_i])\lambda_i)\\
&=& \sum_{i=1}^n Tr_* (\langle {(-1)}^l p_i^{\otimes (2l+1)}
\rangle_{\cal HC})\lambda_i.
\end{array}
$$

On the other hand, we can suppose that
${\cal E}^{(k)}=\{\{\lambda_1\},\dots ,\{\lambda_n\}\}$
and $\tilde a_k=\sum_{i=1}^n p_i^{\otimes (2l+1)}\lambda_i$
for all $k\ge 1$.
Whence,
$$
\begin{array}{ccl}
{\cal C}h_0^{2l}([a])={(-1)}^l T(a)&=&{(-1)}^l
\lim_k Tr_*\langle\tilde a_k\rangle_{\cal HC}\\
&=&{(-1)}^l\sum_{i=1}^n Tr_* (\langle p_i^{\otimes (2l+1)}
\rangle_{\cal HC})\lambda_i.
\pcom
\end{array}
$$

In particular, Theorem~\ref{th8} implies that
one can extend the generalized Chern character to the map from
the quotient group $N_0(A)/\Ker h$ to the even Banach homology.

\section{Generalized Lefshetz numbers}
Suppose $A$, as above, is a von Neumann algebra,
$G$ is a compact Lie group, and $X$ is a compact $G$-manifold.
Let us denote by ${\cal P}(A)$ the category of finitely
generated projective modules over $A$.

Let us recall some notation from~\cite{Sol-Tro}.
The set of all $G$-$A$-bundles over $X$ is an abelian semigroup
with respect to the direct sum operation. The symmetrization
of this semigroup is denoted by $K_G(X;A)$. Assume
$K^G(A):=K_G(pt;A)$.
In this situation there is an isomorphism
\begin{equation}\label{eq11}
K^G(A)\cong K_0(A)\otimes R(G),
\end{equation}
where $R(G)$ is the ring of representations for $G$.

Let us consider a sequence $\{E^i\}$ of $G$-$A$-bundles
over $X$ together with equivariant
pseudo-differential operators
$\{d_i:\Gamma (E^i)\lra \Gamma (E^{i-1})\}$,
where by $\Gamma (E^i)$ we denote
the Banach (with respect to the uniform topology) $A$-module
of continuous sections of $E^i$.
Besides, let us denote by $\sigma _i$ the symbol of $d_i$.
Then this sequence of bundles and operators is called
a {\em $G$-$A$-elliptic complex\/} (and is denoted by $(E,d)$)
if it satisfies the following conditions:\\

(i) $d_i d_{i+1}=0,$\\

(ii) the sequence of symbols
$$
0\lra\pi ^*E^n\stackrel{\sigma _n}{\lra}
\pi ^*E^{n-1}\lra\dots\stackrel{\sigma _1}{\lra}
\pi ^*E^0\lra 0
$$
is exact out of some compact neighbourhood of
the zero section $X\subset T^*X$. Here
$\pi :T^*X\lra X$ is the natural projection.

The index of the elliptic operator
$F=d+d^*:\Gamma (E_{ev})\lra\Gamma (E_{od})$
is an element of the group $K^G(A)$.
Furthermore, for any $g\in G$ by computation
of the character we can define the map
$g:R(G)\lra {\bf C}$. Whence, using isomorphism~(\ref{eq11}),
we obtain the map
$$g:K^G(A)\lra K_0(A)\otimes {\bf C}.$$

Then the {\em Lefschetz number of the first type\/} is
defined as follows
$$L_1(E,g)=g({\rm index}(F))\in K_0(A)\otimes {\bf C}.$$

Note that there exists a connection between these Lefschetz numbers
and fixed points of $g$ (see~\cite{Tro,Tro3}).

Now let us consider an $A$-elliptic complex $(E,d)$
and its unitary endomorphism $U$.
Furthermore, we shall assume that $U=U_g$ for some
representation $U_g$ of a compact Lie group $G$.

Let ${\cal M}$ be a Hilbert $A$-module (see, for example,~\cite{Pas}).
We denote by ${\rm End}_A({\cal M})$ the Banach
algebra of all bounded $A$-homomorphisms of ${\cal M}$.
Now let us consider some strongly continuous representation
$G\lra {\rm End}_A({\cal M})$.
Then this representation is called {\em unitary,\/} if
$\langle gx,gy\rangle = \langle x,y\rangle$
for any $g\in G, x,y\in {\cal M}$.
A Hilbert $A$-module together with a unitary representation
of the group $G$ is called a Hilbert $G$-$A$-module.
Besides, a set $\{x\}_{\beta\in B}\subset {\cal M}$ is
a {\em system of generators\/} for ${\cal M}$, if
finite sums $\{\sum_k x_ka_k : a_k\in A\}$ are dense in~${\cal M}$.

We need the following result of~\cite{Mish}.
\begin{th}\label{th3}
Let ${\cal M}$ be a countably generated Hilbert $G$-$A$-module.
Besides, let $\{V_\pi\}$ be a full system of pairwise
not isomorphic unitary finite-dimentional irreducible
representations for $G$. Then there exists a
$G$-$A$-isomorphism
$$
{\cal M}\cong\bigoplus_{\pi} {\rm Hom}_G(V_\pi,{\cal M})
\otimes_{\bf C} V_\pi.
$$
Here the algebra $A$ (the group $G$) acts on the first
(on the second) multiplier of the
space ${\rm Hom}_G(V_\pi,{\cal M})\otimes_{\bf C} V_\pi$.
$\pcom$
\end{th}

Let the $G$-$A$-module ${\cal M}$ belong to the class
${\cal P}(A)$. Then it is clear that
${\cal M}_\pi:={\rm Hom}_G(V_\pi,{\cal M})\in {\cal P}(A)$.
Furthermore, only finite number of terms in the sum
$\bigoplus\limits_{\pi} {\cal M}_\pi\otimes_{\bf C} V_\pi$
is not equal to zero (see~\cite[1.3.49]{Sol-Tro}).

In particular, we obtain for ${\cal M}=A^n$
the following formula
$$
A^n\cong\bigoplus_{k=1}^MQ_k\otimes V_k,
$$
where $V_k\cong{\bf C}^{L_k}$ and $Q_k\in {\cal P}(A)$.
Therefore,
\begin{equation}\label{eq13}
U_g(\sum_{k=1}^M x_k\otimes v_k)=
\sum_{k=1}^M x_k\otimes u_g^kv_k=
\sum_{k=1}^M\sum_{s=1}^{L_k}x_k\otimes
e^{i\varphi_s^k} v_k^sf_s.
\end{equation}
Here $f_1,\dots ,f_{L_k}$ is a basis for $V_k$
such that the operator $u_g^k$ is diagonal
with respect to it;
$v_k=\sum v_k^sf_s$.
In this case let us define
$$\tau (U_g)=\sum_{k=1}^M
Ch^0_{2l}[Q_k]\cdot {\rm Trace}(u_g^k)\in HC_{2l}(A).
$$

The following result was proved in~\cite{Fra-Tro}.
\begin{lem}
For the $A$-Fredholm operator
$F=d+d^*:\Gamma (E_{ev})\lra\Gamma (E_{od})$
there exists a decomposition
\begin{equation}\label{eq12}
F: M_0\oplus
\widetilde {N_0}\lra M_1\oplus \widetilde {N_1},\; F:M_0\cong M_1
\end{equation}
such that
$$ \widetilde {N_0}=\bigoplus_{j=0}^T N_{2j},\;
\widetilde {N_1}=\bigoplus_{j=1}^T N_{2j-1},\;
N_m\subset\Gamma (E_m),$$
where $N_m$ are projective $U$-invariant Hilbert $A$-modules.
$\pcom$
\end{lem}

Now the {\em Lefschetz number of the second type\/} is
defined as follows
$$L_{2l}(E,U_g)=\sum_j {(-1)}^j\tau (U_g|N_j)\in HC_{2l}(A).$$
This definition is well.

For more detail about W*-Lefschetz numbers we refer
to the works~\cite{Fra-Tro, Tro1, Tro2, Tro3}.

Now let us consider an $A$-elliptic complex $(E,d)$ and an
arbitrary unitary endomorphism $U$ of it
($U$ is not necessarily an element of some representation of $G$).
In this situation let us formulate the following
\begin{dfn}
{\rm We define the} generalized Lefschetz number
 ${\cal L}_1$ {\rm as follows}
$${\cal L}_1(E,U)=\sum_j {(-1)}^j[U|N_j]\in N_0(A).$$
\end{dfn}

Note that generalized Lefschetz numbers are well defined.
This follows by the same reason that for the Lefschetz
numbers of the second type (see~\cite[5.2.21]{Sol-Tro}).

\begin{th}\label{th4}
Let $U$ be a unitary endomorphism of an $A$-elliptic complex
$(E,d)$. Besides, suppose that $U=U_g$ for some
representation $U_g$ of a compact Lie group $G$.
Then ${\cal L}_1(E,U_g)$ belongs to the group $N_0(A)_{fin}$
and $h({\cal L}_1(E,U_g))=L_1(E,U_g)$,
where $h$ is map~{\rm (\ref{eq10})}.
\end{th}

{\bf Proof.}
Let us examine decomposition~(\ref{eq12}) for $F$.
We have shown above that there exist isomorphisms
$$N_{2j}\cong\bigoplus_{k_j=1}^{K_j} P_{k_j}^{(j)}\otimes
V_{k_j}^{(j)},\quad
N_{2j-1}\cong\bigoplus_{l_j=1}^{L_j} Q_{l_j}^{(j)}\otimes
W_{l_j}^{(j)},
$$
where $V_{k_j}^{(j)}$ and $W_{l_j}^{(j)}$ are
comlex vector spaces of irreducible unitary representations
of $G$, $P_{k_j}^{(j)}$ and $Q_{l_j}^{(j)}$ are
$G$-trivial modules from ${\cal P}(A)$.
Thus we get
$${\rm index}(F)=\sum_{j=0}^T\sum_{k_j=1}^{K_j}
[P_{k_j}^{(j)}]\otimes \chi (V_{k_j}^{(j)})-
\sum_{j=1}^T\sum_{l_j=1}^{L_j}
[Q_{l_j}^{(j)}]\otimes \chi(W_{l_j}^{(j)})
$$
and
$$
L_1(E,g)=\sum_{j=0}^T\sum_{k_j=1}^{K_j}
[P_{k_j}^{(j)}]\otimes {\rm Trace}(g|V_{k_j}^{(j)})-
\sum_{j=1}^T\sum_{l_j=1}^{L_j}
[Q_{l_j}^{(j)}]\otimes {\rm Trace}(g|W_{l_j}^{(j)}).
$$
Here we have denoted by $\chi$ the character of the
representation of $G$.

On the other hand, using expression~(\ref{eq13}), we conclude
that
$$sp(U_g|(P_{k_j}^{(j)}\otimes V_{k_j}^{(j)}))=
sp(u_g^{(j),k_j}|V_{k_j}^{(j)}).$$
This implies that ${\cal L}_1(E,U_g)$ belongs to the group
$N_0(A)_{fin}$.
Furtermore, for any $e^{i\varphi}$ from
$sp(u_g^{(j),k_j}|V_{k_j}^{(j)})$ the spectral projection
of $U_g$ corresponding to this point is equal to
$P_{k_j}^{(j)}$.
Thus we obtain
$$
\begin{array}{ccl}
h({\cal L}_1(E,U_g))&=&h(\sum_{j=0}^{2T} {(-1)}^j
[U_g|N_j])\\
&=& h(\sum_{j=0}^T\sum_{k_j=1}^{K_j}
[Id_{P_{k_j}^{(j)}}\otimes u_g^{(j),k_j}|V_{k_j}^{(j)}]\\
&-&\sum_{j=1}^T\sum_{l_j=1}^{L_j}
[Id_{Q_{l_j}^{(j)}}\otimes u_g^{(j),l_j}|W_{l_j}^{(j)}])\\
&=&h(\sum_{j=0}^T\sum_{k_j=1}^{K_j}
[P_{k_j}^{(j)}{\rm Trace} (u_g^{(j),k_j}|V_{k_j}^{(j)})]\\
&-&\sum_{j=1}^T\sum_{l_j=1}^{L_j}
[Q_{l_j}^{(j)}{\rm Trace}(u_g^{(j),l_j}|W_{l_j}^{(j)})])\\
&=& \sum_{j=0}^T\sum_{k_j=1}^{K_j}
[P_{k_j}^{(j)}]\otimes {\rm Trace}(u_g^{(j),k_j}|V_{k_j}^{(j)})\\
&-&\sum_{j=1}^T\sum_{l_j=1}^{L_j}
[Q_{l_j}^{(j)}]\otimes {\rm Trace}(u_g^{(j),l_j}|W_{l_j}^{(j)}).
\end{array}
$$
Thus we establish the required statement.
$\pcom$

\begin{th}\label{th5}
Under the assumptions of the previous theorem,
we have
$$\pi_*(L_{2l}(E,U_g))={\cal C}h^0_{2l}({\cal L}_1(E,U_g)).$$
Here $\pi_*$ is map {\rm (\ref{eq2})}.
\end{th}

{\bf Proof.}
From Theorems~\ref{th8},~\ref{th4} we deduce that
$$
\begin{array}{ccl}
{\cal C}h^0_{2l}({\cal L}_1(E,U_g)) &=&
\pi_*\,\widetilde{Ch}_{2l}^0\, h({\cal L}_1(E,U_g))\\
&=& \pi_*\widetilde {Ch}_{2l}^0 (L_1(E,g)).
\end{array}
$$
Furtermore, it follows from~\cite[Theorem 5.2.22]{Sol-Tro}
that $\widetilde {Ch}_{2l}^0 (L_1(E,g))=L_{2l}(E,U_g)$.
The proof is complete.
$\pcom$ \\

{\bf Acknowledgements.} The author is grateful to Prof.
E.V. Troitsky for the attention to the work and also to
Dr. V.M. Manuilov and Prof. A.S.~Mishchenko
for helpful discussions.


 \begin{center}
Alexandre Pavlov\\
Chair of Higher Geometry and Topology\\
Dept. of Mech. and Mathematics\\
Moscow State University\\
Vorobjevi gori\\
Moscow 119899\\
Russia\\
email: pavlov@mech.math.msu.su
 \end{center}

\end{document}